\documentclass{article}
\usepackage{amssymb}
\title{ON THE SZ\"USZ'S SOLUTION TO \\ GAUSS' PROBLEM}
\author{
	Ion Col\c tescu, Dan Lascu\\
	\small ``Mircea cel Batran'' Naval Academy, 1 Fulgerului,\\
	\small 900218 Constanta, Romania\\
	\small {\it e-mail}: icoltescu@yahoo.com, lascudan@gmail.com 
}
\date{}
\newcommand{\N}{\mathbb N}
\newcommand{\Q}{\mathbb Q}
\newcommand{\C}{\mathbb C}
\begin{document}
\maketitle
\begin{abstract}
The present paper deals with Gauss' problem on continued fractions. We present a new proof of a theorem which Sz\"usz applied in order to solve this problem.
To be noted, that we obtain the value $0.7594\dots$ for $q$, which has been optimized by Sz\"usz in his 1961 paper "\"Uber einen Kusminschen Satz", where the value $0.485$ is obtained for $q$.
In our proof, we make use of an important property of the Perron-Frobenius operator of $\tau$ under $\gamma$, where $\tau$ is the continued fraction transformation, and $\gamma$ is the Gauss' measure.

\end{abstract}

\section{INTRODUCTION}

Let $\xi \in [0,1)$, and let
\[
\xi = \frac{1}{\displaystyle d_1 + \frac{1}{\displaystyle d_2 + \ddots + \frac{1}{\displaystyle d_n+\ddots}}}=[0;d_1,d_2,\ldots,d_n,\ldots]
\]
be the regular continued fraction expansion of $\xi$.
On October 25, 1800, Gauss wrote in his diary that (in modern notation):
\begin{equation}
\lim_{n\rightarrow\infty}\lambda\left(\left\{\xi\in[0,1);\tau^n(\xi) \leq z \right\}\right)=\frac{\log(z+1)}{\log 2}, 0\leq z\leq 1,
\end{equation}
where $\lambda$ is the Lebesgue measure and $\tau:[0,1)\rightarrow[0,1)$ is the continued fraction transformation defined by  
\begin{equation}
\tau(\xi):=\frac{1}{\xi}-\left[\frac{1}{\xi}\right], \xi \neq 0; \tau(0):=0,
\end{equation}
where $[\cdot]$ denotes the floor (entire) function.
Latter on, in a letter dated January 30, 1812, Gauss asked Laplace to provide an estimate of the error term $r_n(z)$, defined by  
\[
r_n(z) := \lambda(\tau^{-n}([0,z]))-\frac{\log(z+1)}{\log 2}, n \geq 1.
\]

Gauss' proof has never been found.
The first who prove (1) and at the same time answered to Gauss' question was Kuzmin.
In 1928, Kuzmin [3] showed that  
\[
r_n(z)={\cal O}(q^{\sqrt n}),
\]
with $q\in(0,1)$, uniformly in $z$.
Independently, Paul L\'evy showed one year later that  
\[
r_n(z)={\cal O}(q^n),
\]
with $q=0.7\ldots$, uniformly in $z$.
From that moment onwards, a great number of such Gauss-Kuzmin theorems followed.
To mention a few: F.Schweiger (1968), P. Wirsing [6] (1974 - which determined that the optimal value of $q$ is equal to $0.303663002$), K.I. Babenko (1978), and more recently M. Iosifescu (1992).

\section{THE GAUSS-KUZMIN TYPE EQUATION}
An essential ingredient in any proof of whichever Gauss-Kuzmin theorem is the following observation.
Let $\xi\in[0,1)\backslash\Q$ and put
\[
    \tau_k := \tau^k(\xi), k\geq 0,
\]
where $\tau:[0,1)\rightarrow[0,1)$ is the continued fraction transformation defined in (2).
From (2) it follows that  
\[
    0\leq \tau_{n+1} \leq x \Leftrightarrow \tau_n \in \bigcup_{i\in \N_+}\left[\frac{1}{x+i},\frac{1}{i}\right].
\]
Thus, if we put
\[
    F_n(x):=\lambda\left(\left\{\xi\in[0,1);\tau^n(\xi)\leq x\right\}\right), n\geq 0,
\]
then  
\begin{equation}
F_{n+1}(x)=\sum_{i\in\N_+}\left(F_n\left(\frac{1}{i}\right)-F_n\left(\frac{1}{x+i}\right)\right), n\geq 0,
\end{equation}
relation called the Gauss-Kuzmin type equation.

\section{IMPORTANT RESULT}
Let $B(I)$ the Banach space of all bounded measurable functions $f:I\rightarrow\C$, $I:=[0,1]$.

\noindent {\bf Proposition}. If $f\in B(I)$ is non-decreasing, then $Uf$ is non-increasing, where $U$ is the Perron-Frobenius operator of $\tau$ under $\gamma$, with $\gamma$ the Gauss' measure which is defined on ${\cal B}_{[0,1]}$ - Borel $\sigma$-algebra of sets on $[0,1]$, by 
\[
\gamma(A)=\frac{1}{\log 2}\int_A\frac{dx}{x+1}, A\in{\cal B}_{[0,1]}.
\]

\noindent {\bf Proof}. Let $f$ be a non-decreasing function. Thus, if $x<y$, then $f(x)\leq f(y)$.
We evaluate the difference $Uf(y)-Uf(x)$.
We have, $\displaystyle Uf(y)=\sum_{i\in \N_+}P_i(y)f\left(\frac{1}{y+i}\right)$ and $\displaystyle Uf(x)=\sum_{i\in \N_+}P_i(x)f\left(\frac{1}{x+i}\right)$, where $\displaystyle P_i(x) = \frac{x+1}{(x+i)(x+i+1)}$.
Thus, $Uf(y)-Uf(x) = S_1 + S_2$, where 
\[
  S_1=\sum_{i\in \N_+}P_i(y)\left(f\left(\frac{1}{y+i}\right) - f\left(\frac{1}{x+i}\right)\right),
  S_2=\sum_{i\in \N_+}(P_i(y)-P_i(x))f\left(\frac{1}{x+i}\right).
\]
Since $f$ is non-decreasing, and $\frac{1}{x+i}>\frac{1}{y+i}$, then $f\left(\frac{1}{x+i}\right)\geq f\left(\frac{1}{y+i}\right)$.\\
Thus, $S_1\leq 0$.
We will show that $S_2\geq 0$ too. 
We have that  $\displaystyle \sum_{i\in\N_+}P_i(u)=1$, $u\in I$, and therefore we obtain: 
\[
\begin{array}{rl}
S_2 = &\displaystyle \sum_{i\in\N_+}(P_i(y)-P_i(x))f\left(\frac{1}{x+i}\right) - \sum_{i\in\N_+}(P_i(y)-P_i(x))f\left(\frac{1}{x+1}\right)\\
=&\displaystyle 
-\sum_{i\in\N_+}\left(f\left(\frac{1}{x+1}\right)-f\left(\frac{1}{x+i}\right)\right)(P_i(y)-P_i(x)).
\end{array}
\]
Now, it is easy to show that the function $P_1$ is decreasing, while the functions $P_i$, $i\geq 2$, are all increasing. 
Also, 
\[
f\left(\frac{1}{x+1}\right)-f\left(\frac{1}{x+i}\right)\geq f\left(\frac{1}{x+1}\right)-f\left(\frac{1}{x+2}\right)\geq 0, i\geq 2.
\]
Therefore,
\[
\begin{array}{rl}
S_2 = &\displaystyle -\sum_{i\geq 2}\left(f\left(\frac{1}{x+1}\right)-f\left(\frac{1}{x+i}\right)\right)(P_i(y)-P_i(x)) \\
\leq  &\displaystyle 
- \left(f\left(\frac{1}{x+1}\right)-f\left(\frac{1}{x+2}\right)\right)\sum_{i\geq 2}(P_i(y)-P_i(x)) \leq 0.
\end{array}
\] 
Thus, $Uf(y)-Uf(x)\leq 0$.

\section{THE GAUSS-KUZMIN THEOREM}
We will give a simple proof that 
\[
F_n(x)=\frac{\log(x+1)}{\log 2}+{\cal O}(q^n),
\]
where $0<q<1$ or, to be exactly, $q=0.7594\ldots$. 
In fact, we will proof the following:

\noindent {\bf Theorem}. Let $f_0(x)$ be any twice differentiable function defined on $[0,1]$ with $f_0(0)=0$ and $f_0(1)=1$.
Let the sequence of functions $f_1(x), f_2(x), \ldots$ be defined by the recursion formula 
\[
f_{n+1}(x) = \sum_{i\in\N_+}\left(f_n\left(\frac{1}{i}\right)-f_n\left(\frac{1}{x+i}\right)\right).
\]
Then 
\[
f_{n}(x) = \frac{\log(x+1)}{\log 2} + {\cal O}(q^n),
\]
where $0<q<1$ or, to be exactly, $q=0.7594\ldots$.

It is clear that for $f_0(x)=x=F_0(x)$, this theorem will establish Gauss' claim and provide an answer to his problem.

\noindent {\bf Proof}. Instead of studying $f_n(x)$ directly, we look at the derivative:
\begin{equation}
	f'_{n+1}(x) = \sum_{i\in\N_+}\frac{1}{(x+i)^2}f'_n\left(\frac{1}{x+i}\right).
\end{equation}

Let us introduce another sequence of functions $g_0,g_1,\ldots$ defined by 
\[
	g_n(x)=(x+1)f'_n(x).
\]
Then the recursion formula (4) is transformed into 
\[
	\frac{g_{n+1}(x)}{x+1} = \sum_{i\in\N_+}\frac{1}{(x+i)^2}\frac{g_n\left(\frac{1}{x+i}\right)}{\frac{1}{x+i}+1} 
	= \sum_{i\in\N_+}\frac{1}{(x+i)(x+i+1)}g_n\left(\frac{1}{x+i}\right) \Rightarrow
\]
\[
\Rightarrow g_{n+1}(x) = \sum_{i\in\N_+}\frac{x+1}{(x+i)(x+i+1)}g_n\left(\frac{1}{x+i}\right) 
	= \sum_{i\in\N_+}P_i(x)g_n\left(\frac{1}{x+i}\right) = Ug_n,
\]
where $P_i(x) = \frac{x+1}{(x+i)(x+i+1)}$, $i\in\N_+$, and $U$ is the Perron-Frobenius operator of $\tau$ under $\gamma$. 

If we can show that $g_n(x)=\frac{1}{\log 2}+{\cal O}(g^n)$, then an integration will establish the theorem for $f_n(x)$, because integrating $\frac{1}{x+1}$ will give $\log(x+1)$ term together with a bounded expression on a bounded interval times the ${\cal O}(q^n)$ error term, which will remain ${\cal O}(q^n)$.
To demonstrate that $g_n(x)$ has this desired form, it suffices to establish that $g'_n(x) = {\cal O}(q^n)$, as the $\frac{1}{\log 2}$ constant in $g_n(x)$ will follow from the normalization requirement that $f_0(0)=0$ and $f_0(1)=1$.

We have: 
\[
	P_i(x) = \frac{x+1}{(x+i)(x+i+1)} = \frac{i}{x+i+1} - \frac{i-1}{x+i},
\]
thus 
\[
	g_{n+1}(x) = \sum_{i\in\N_+} \left(\frac{i}{x+i+1} - \frac{i-1}{x+i}\right)g_n\left(\frac{1}{x+i}\right) \Leftrightarrow
\]
\begin{equation}
\begin{array}{rl}
g'_{n+1}(x) =&\displaystyle - \sum_{i\in\N_+} \left(\frac{i}{(x+i+1)^2}\right) \left(g_n\left(\frac{1}{x+i}\right)-g_n\left(\frac{1}{x+i+1}\right)\right)- \\
&\displaystyle - \sum_{i\in\N_+} \frac{P_i(x)}{(x+i)^2} g'_n\left(\frac{1}{x+i}\right).
\end{array}
\end{equation}
By applying the mean value theorem of calculus to the difference 
\[
g_n\left(\frac{1}{x+i}\right) - g_n\left(\frac{1}{x+i+1}\right),
\]
we obtain 
\[
	g_n\left(\frac{1}{x+i}\right) - g_n\left(\frac{1}{x+i+1}\right) = \left(\frac{1}{x+i} -\frac{1}{x+i+1} \right) g'_n\left(\frac{1}{x+\theta_i}\right),
\]
where $1<\theta_i<i$.

Thus, from (5), we have:
\begin{equation}
\begin{array}{rl}
	g'_{n+1}(x) =&\displaystyle -\sum_{i\in\N_+} \frac{i}{(x+i)(x+i+1)^3} g'_n\left(\frac{1}{x+\theta_i}\right) -\\
	&\displaystyle -\sum_{i\in\N_+} \frac{P_i(x)}{(x+i)^2} g'_n\left(\frac{1}{x+i}\right)
\end{array}
\end{equation}
Let $M_n$ be the maximum of $|g'_n(x)|$ on $[0,1]$, i.e. $\displaystyle M_n=\max_{x\in[0,1]}|g'_n(x)|$.\\
Then, from (6), we have that:
\begin{equation}
	M_{n+1} \leq M_n \max_{x\in[0,1]}\left(
	\sum_{i\in\N_+} \frac{i}{(x+i)(x+i+1)^3} + \sum_{i\in\N_+} \frac{P_i(x)}{(x+i)^2}
	\right).
\end{equation}
We now must calculate the maximum value of the sums in this expression. 
To this end, define function $h$ by 
\[
	h(x) = \sum_{i\in\N_+} \frac{P_i(x)}{(x+i)^2}, x\in [0,1].
\]
Note that for $\varphi(x)=x^2$, $x\in[0,1]$, we have $h(x)=U\varphi(x)$.
Since $\varphi$ is increasing, and using the proposition from Section 3, we have that $h$ is decreasing. 
Hence, $h(x)\leq h(0)$, and since $\displaystyle h(0) = \sum_{i\in\N_+} \frac{P_i(0)}{i^2} = \sum_{i\in\N_+} \frac{1}{i^3 (i+1)}$.
Therefore, (7) become: 
\[
\begin{array}{rl}
	M_{n+1} \leq & \displaystyle M_n \sum_{i\in\N_+} \left(\frac{1}{i(i+1)^3} + \frac{1}{i^3(i+1)}\right) \\
	
	=&\displaystyle M_n \sum_{i\in\N_+} \left(\frac{1}{i^3} - \frac{1}{i^2} + \frac{1}{i} - \frac{1}{i+1}+\frac{1}{(i+1)^3}\right) \\
	=&\displaystyle  M_n (\zeta(3) - \zeta(2) + 1 + \zeta(3) - 1) \\
	=&\displaystyle  M_n (2\zeta(3) - \zeta(2)),
\end{array}
\]
where $\zeta(n)$ denotes the Riemann zeta function.
Hence, $2\zeta(3)-\zeta(2)=0.7594798\ldots$.
Thus, $M_{n+1} < M_n q$, where $q=0.7594798\ldots$, and $q'_{n+1}(x) = {\cal O}(q^{n+1})$, which proves the theorem.

\end{document}